\documentclass[conference]{IEEEtran}
\IEEEoverridecommandlockouts
\usepackage{cite}
\usepackage{amsmath,amssymb,amsfonts}
\usepackage{algorithmic}
\usepackage{graphicx}
\usepackage{textcomp}
\usepackage{xcolor}
\usepackage[numbers,sort&compress]{natbib}
\usepackage{booktabs}
\PassOptionsToPackage{hyphens}{url}
\usepackage[hidelinks]{hyperref}
\usepackage{multirow}

\def\BibTeX{{\rm B\kern-.05em{\sc i\kern-.025em b}\kern-.08em
    T\kern-.1667em\lower.7ex\hbox{E}\kern-.125emX}}
\begin{document}

\title{Energy-Efficient Routing for Electric Vehicles under Acceleration and Load Effects}

\author{\IEEEauthorblockN{Tingting Su}
\IEEEauthorblockA{\textit{School of Management} \\
\textit{Xiamen University} \\
Xiamen, China \\
sutingting@stu.xmu.edu.cn}
\and
\IEEEauthorblockN{Xinyue Zhang}
\IEEEauthorblockA{\textit{School of Computer Science and Engineering} \\
\textit{Beihang University}\\
Beijing, China\\
zhangxinyue.hihi@gmail.com}
\and
\IEEEauthorblockN{Jingyi Zhao*}
\IEEEauthorblockA{\textit{Shenzhen Research Institute of Big Data} \\
Shenzhen, China\\
zhaojingyi@sribd.cn}
}
\IEEEaftertitletext{\vspace{-1.0\baselineskip}}
\maketitle
\begin{abstract}
This paper proposes an Acceleration and Load-Dependent Electric Vehicle Routing Problem (ALD-EVRP), to optimize the energy consumption (EC) while capturing the effects of changing traffic conditions between peak and off-peak periods. We generalize the time-dependent speed model by replacing step functions with piecewise linear functions. The EC of each vehicle is influenced by its speed, acceleration, and real-time load. A mathematical model is developed and solved using BonMin, and a custom meta-heuristic algorithm is proposed for large-scale problems, yielding the same results as BonMin on small problems and performing better on larger ones. This is validated with real data from Singapore.
\end{abstract}


\section{Introduction}
In recent years, the number of Electric Vehicles (EVs) has increased rapidly, from a few thousand in 2009 to about $6.8$ million in use globally by the end of 2020. This surge in EV adoption has spurred increased interest in the Electric Vehicle Routing Problem (EVRP), leading to the development of numerous EVRP variants under diverse practical constraints, objectives, and decisions.
In modern urban areas, vehicle speeds vary significantly on individual roads between peak and off-peak periods. Therefore, it is essential to include time-dependent speeds in routing plans. Real-time road conditions limit vehicle speeds, leading to the incorporation of time-dependent conditions in the traditional vehicle routing problem (VRP) to make research more applicable to real-world settings. Previous research frequently models speed using a stepwise speed function, implying sudden changes at specific points, which fails to capture the gradual variations in vehicle speed that occur when traffic shifts from peak to off-peak periods.
To tackle this problem, we develop a time-dependent (TD) speed model. In this model, vehicle speed on each arc changes gradually, represented by a {\em piecewise linear function}, reflecting real-time road conditions. Furthermore, the energy consumption (EC) of EVs, crucial in the EVRP, varies with the vehicle's load, affecting both route planning and vehicle efficiency \citep{KANCHARLA2020113714,BEHNKE2021794}.
Motivated by \citep{naeem2024energy}, which showed that jointly optimizing route selection and driving speed strategies can substantially reduce energy consumption and improve battery performance, we incorporate a load-dependent term into our energy-consumption model to more accurately reflect the effect of payload on total energy use and enhance the realism of route and energy evaluations.

Considering LD factors enhances the EVRP by enabling a more accurate simulation of EV energy consumption. This improvement allows for better route planning, recharging strategies, and cargo distribution, which together minimize energy costs, extend battery life, and reduce maintenance expenses. Our proposed Acceleration and Load-Dependent EVRP (ALD-EVRP) accounts for the impact of real-time load and acceleration on EC. The goal is to reduce the overall EC, which depends on the vehicle’s carried load, instantaneous speed, and acceleration. For simplicity, this model excludes the recharging process.
The contributions of this paper are as follows.
\begin{enumerate}
    \item We formulate a new problem {(ALD-EVRP)} and provide its mathematical model. This model is a Mixed Integer Nonlinear Programming (MINLP) framework that incorporates real-time road conditions and load effects. The accuracy of the model is validated using the BonMin Optimizer.
    \item We develop a meta-heuristic algorithm combining Large Neighborhood Search (LNS) and Local Search (LS), using the split algorithm \citep{PRINS20041985} to partition fleet tours into individual EV routes.
    \item We generate 47 ALD-EVRP test instances for up to 114 customers using data from Singapore, and evaluate our LNS-SPP algorithm. 
\end{enumerate}

The remainder of the paper is organized as follows. First, the problem is introduced along with its mathematical formulation. Next, a novel meta-heuristic algorithm, specifically designed to address the challenges of this problem, is presented. Finally, the experimental data are reported, the algorithm's performance is evaluated on both large-scale and small-scale test cases, and the results are discussed.

\section{Mathematical Formulation}
The ALD-EVRP with a fleet of vehicles is modeled on a directed graph $\mathcal{G} = (\mathcal{V}, \mathcal{B})$ with nodes $\mathcal{V} = \{0,1,\dots,n+1\}$, where $0$ and $n+1$ are depots, and $1$ to $n$ are customers with demands $q_i$ set to $0$ for depots. For convenience, we define three subsets: 
$\mathcal{V}^+ = \mathcal{V}^+$, which contains all nodes except the starting depot; 
$\mathcal{V}^- = \mathcal{V}^-$, which contains all nodes except the ending depot; 
and $\tilde{\mathcal{V}} = \mathcal{V} \setminus \{0,n+1\}$, which contains only the customer nodes. The arcs $\mathcal{B} = \{(i,j)| i \in \mathcal{V}^{-}, j \in \mathcal{V}^{+}\}$ link nodes with distances $d_{ij}$ and variable speed profiles. A set $\mathcal{K}$ of EVs, each assigned a load capacity $Q$, operates on this network to deliver goods efficiently under varying traffic conditions.

\begin{figure}[htbp]
\centering
\includegraphics[width=0.48\textwidth]{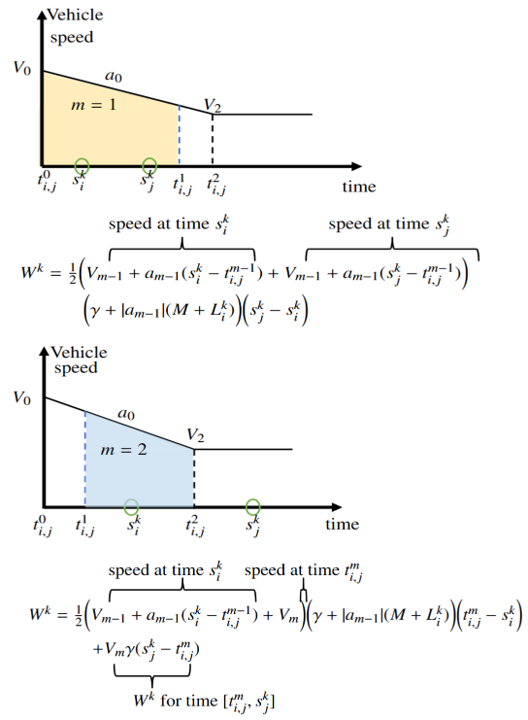}
\caption{Calculation of $W^k$ for $m\in\Phi_1,\Phi_2$ (acceleration and cruising)}
\label{fig:W12}
\end{figure}

\begin{figure}[htbp]
\centering
\includegraphics[width=0.48\textwidth]{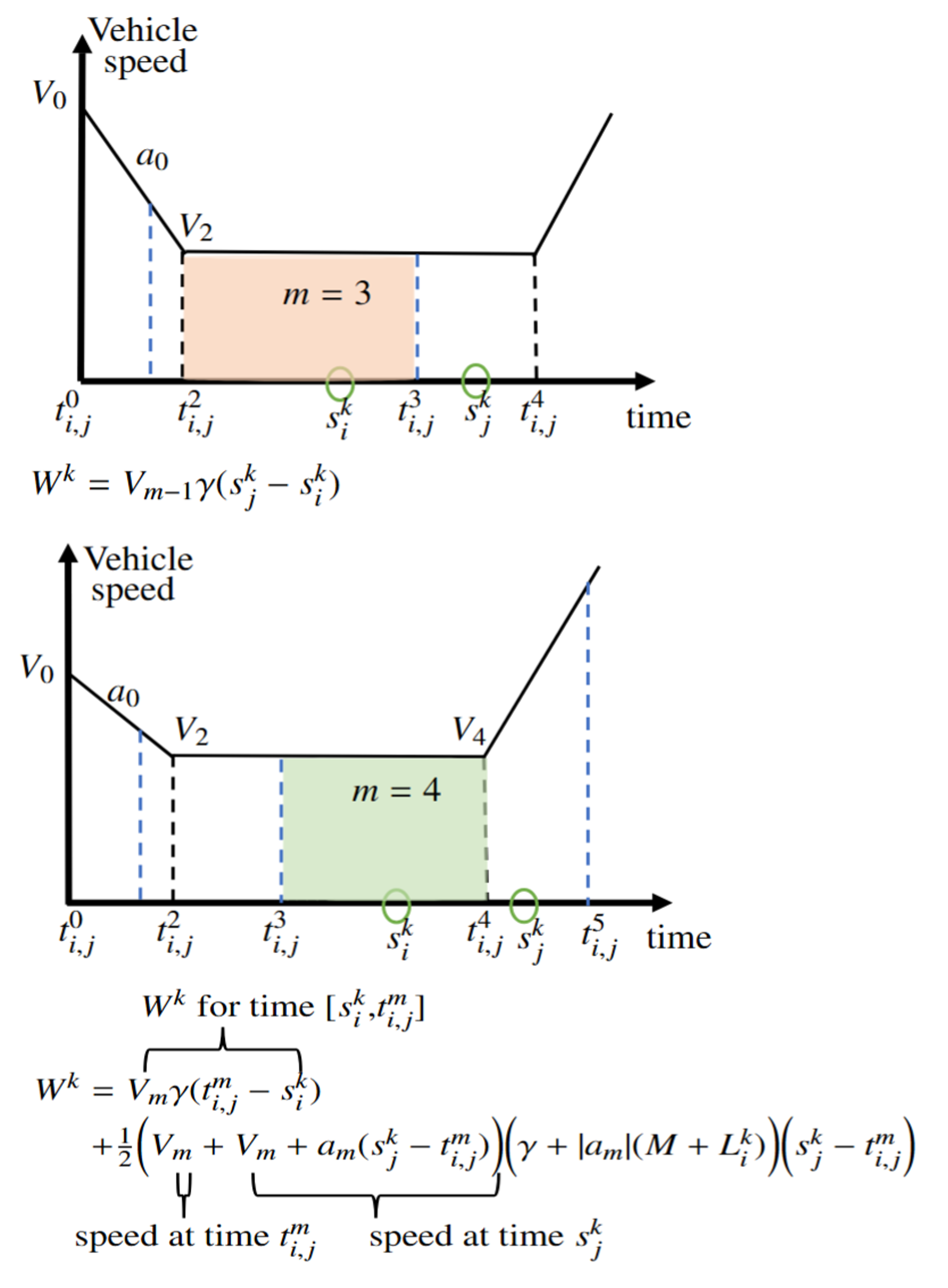}
\caption{Calculation of $W^k$ for $m\in\Phi_3,\Phi_4$ (deceleration and idling)}
\label{fig:W34}
\end{figure}

For each travel arc $(i,j)\in\mathcal{B}$, the driving process is segmented into consecutive time intervals to capture changes in traffic conditions and vehicle speed. 
Let $T_{i,j}=\{1,2,\dots\}$ denote the ordered set of intervals, which are delimited by breakpoints $\{t_{i,j}^{0}, t_{i,j}^{1}, t_{i,j}^{2},\dots\}$.  Within each interval $m\in T_{i,j}$, the vehicle is assumed to move with constant acceleration $a_{i,j}^{m}$, and the instantaneous speed can be expressed as $V_{i,j}(t) = V_{i,j}^{m-1} + a_{i,j}^m (t - t_{i,j}^{m-1}),  t\in[t_{i,j}^{m-1},t_{i,j}^{m}]$. This piecewise-linear representation captures the non-uniform evolution of speed along $(i,j)$, enabling a realistic representation of traffic variations and a consistent time basis for energy computation. Each interval $m\in T_{i,j}$ reflect different phases of a journey, acceleration, maintaining a constant speed, and deceleration. This classification allows the model to represent different driving behaviors in the energy computation, as shown in Figures~\ref{fig:W12} and~\ref{fig:W34}.

The decision variables used in the model are defined as follows. $x_{i,j}^{m,k}$ is a binary variable for each pair$(i,j) \in \mathcal{B}$ and each time interval$m \in T_{i,j}$, which equals 1 if vehicle $k$ departs from node $i$ to node $j$ during that time interval, and 0 otherwise. The continuous variable $w_{i,j}^{m,k}$  denotes the actual departure time of vehicle  $k$ within that time $T_{i,j}^m$. $y_{i,j}^k$ is another binary variable that equals 1 if vehicle $k$ makes the journey from node $i$ to node $j$ . The decision variable $s_i^k$ denotes the departure time from node $i$  for vehicle $k$, with the initial departure time from the depot (node 0) set to $s_0^k = 0$. The set $\mathcal{V}$ represents the collection of nodes, which includes both the depot and customer locations. The set $\mathcal{B}$ defines the arcs between nodes, and $T_{i,j}$ specifies the set of time intervals for travel between nodes $i$ and $j$ . The set  $\mathcal{K}$ defines the vehicles available for routing. In addition, $L_{i,j}^k$  represents the load of vehicle  $k$  after leaving node $i$, while $W_i^k$  measures the total energy consumed by vehicle $k$  subsequent to serving node $i$ . The model also incorporates constraints related to vehicle capacities: $Q_e$ , the maximum load capacity of a vehicle, and $Q_b$, the maximum energy capacity of a vehicle. These variables and sets collectively enable precise tracking and optimization of vehicle routing and scheduling within the framework of the ALD-EVRP.

We model the energy consumption on arc $(i,j)$ using a power–speed–acceleration relation. 
Building on the formulation in \cite{GALVIN2017234}, we make several simplifying assumptions and account for the influence of real-time load $L$ on energy consumption, leading to the following instantaneous power demand: $P(t) = r V_t + s V_t^2+ c V_t^3+(M + L) |a| V_t.$

Assuming that acceleration $a$ remains constant within the interval $[t_1,t_2]$, the energy consumed over this period is

$W = \int_{t_1}^{t_2} P(t)\,dt = \int_{ t_1}^{ t_2} \big( rV_t + sV_t^2 + cV_t^3 + (M + L)|a|V_t \big) \, \mathrm{d}t \notag$,

where $V_{t_1}$ and $V_{t_2}$ denote the speeds at the start and end of the interval. 
Because $V_t$ varies linearly with time, this integration yields an exact value.  
Thus, each arc $(i,j)$ is divided into subsegments corresponding to the four motion phases $\Phi_1$–$\Phi_4$, and the total energy of the vehicle $k$ is obtained by summing the contributions in all segments.

$W^k = \sum_{(i,j)\in\mathcal{B}}\sum_{m\in T_{i,j}} W_{i,j}^{m,k}$.


The overall optimization objective is to select routes and schedules that minimize the total energy use across all vehicles and arcs.

\begin{align*}
\min &~ \sum_{k\in \mathcal{K}} W_{n+1}^k \label{eq3}
\end{align*}

subject to:
\allowdisplaybreaks
\begin{IEEEeqnarray}{rCl}
\sum_{j \in \mathcal{V}^+} \sum_{m \in T_{0,j}} x^{m,k}_{0,j}
&=&
\sum_{i \in \mathcal{V}^-} \sum_{m \in T_{i,n+1}} x^{m,k}_{i,n+1}, \nonumber \\
\IEEEeqnarraymulticol{3}{c}{\forall\, k \in \mathcal{K}.} \label{eq:4} \\
\nonumber \end{IEEEeqnarray}
\begin{IEEEeqnarray}{rCl}
\sum_{k \in \mathcal{K}} \sum_{(i,r) \in \mathcal{B}} \sum_{m \in T_{i,r}} x^{m,k}_{i,r}
&=&
\sum_{k \in \mathcal{K}} \sum_{(r,j) \in \mathcal{B}} \sum_{m \in T_{r,j}} x^{m,k}_{r,j}
= 1, \nonumber \\
\IEEEeqnarraymulticol{3}{c}{\forall\, r \in \tilde{\mathcal{V}}.} \label{eq:5} \\
 \nonumber \end{IEEEeqnarray}
\begin{equation}\label{eq:11}
y_{i,j}^k = \sum_{m \in T_{i,j}} x_{i,j}^{m,k},\quad \forall\, (i,j) \in \mathcal{B},~ k \in \mathcal{K}.
\end{equation}

\begin{IEEEeqnarray}{rCl}
L_j^k
&\geq&
L_{i,j}^k - q_j + \mathbb{M}\big(y_{i,j}^k - 1\big), \forall\, (i,j) \in \mathcal{B},~ k \in \mathcal{K} \label{eq:12} 
\end{IEEEeqnarray}
\begin{IEEEeqnarray}{rCl}
L_0^k &\leq& Q_e, \forall\, k \in \mathcal{K} \label{eq:capa} \end{IEEEeqnarray}
\begin{IEEEeqnarray}{rCl}
W_{n+1}^k &\leq& Q_b, \forall\, k \in \mathcal{K} 
\label{eq:payload2}
\end{IEEEeqnarray}

\begin{small}
\allowdisplaybreaks

\begin{IEEEeqnarray}{rCl}
W_j^k + \mathbb{M}(1 - x_{i,j}^{m,k})
&\geq&
W_i^k
+ (s_j^k - s_i^k)
\Big[\gamma + |a_{m-1}|(M + L_{i,j}^k)\Big] \nonumber\\
&& {}\times
\Big[V_{m-1} + a_{m-1}\Big(\tfrac{s_i^k + s_j^k}{2} - t_{i,j}^{m-1}\Big)\Big], \nonumber\\
\IEEEeqnarraymulticol{3}{c}{
\forall\, m \in \Phi_1,\,
i \in \mathcal{V}^-,\,
j \in \mathcal{V}^+.}
\label{eq:6_new}
\end{IEEEeqnarray}

\begin{IEEEeqnarray}{rCl}
W_j^k + \mathbb{M}(1 - x_{i,j}^{m,k})
&\geq&
W_i^k
+ (t_{i,j}^{m} - s_i^k)
\Big[\gamma + |a_{m-2}|(M + L_{i,j}^k)\Big] \nonumber\\
&& {}\times
\Big[V_m + a_{m-2}\Big(\tfrac{V_{m-2} - V_m}{2a_{m-2}} 
+ s_i^k - t_{i,j}^{m-2}\Big)\Big] \nonumber\\
&& {}+ V_m \gamma (s_j^k - t_{i,j}^{m}), \nonumber\\
\IEEEeqnarraymulticol{3}{c}{
\forall\, m \in \Phi_2,\,
i \in \mathcal{V}^-,\,
j \in \mathcal{V}^+.}
\label{eq:62_new}
\end{IEEEeqnarray}

\begin{IEEEeqnarray}{rCl}
W_j^k + \mathbb{M}(1 - x_{i,j}^{m,k})
&\geq&
W_i^k
+ \gamma V_{m-1}(s_j^k - s_i^k), \nonumber\\
\IEEEeqnarraymulticol{3}{c}{
\forall\, m \in \Phi_3,\,
i \in \mathcal{V}^-,\,
j \in \mathcal{V}^+.}
\label{eq:63_new}
\end{IEEEeqnarray}

\begin{IEEEeqnarray}{rCl}
W_j^k + \mathbb{M}(1 - x_{i,j}^{m,k})
&\geq&
W_i^k
+ (s_j^k - t_{i,j}^{m})
\Big[\gamma + |a_m|(M + L_{i,j}^k)\Big] \nonumber\\
&& {}\times
\Big[V_m + a_m\Big(\tfrac{s_j^k - t_{i,j}^{m}}{2}\Big)\Big]
+ V_m \gamma (t_{i,j}^{m} - s_i^k), \nonumber\\
\IEEEeqnarraymulticol{3}{c}{
\forall\, m \in \Phi_4,\,
i \in \mathcal{V}^-,\,
j \in \mathcal{V}^+.}
\label{eq:64_new}
\end{IEEEeqnarray}

\begin{IEEEeqnarray}{rCl}
s_j^k 
&\geq&
s_i^k + \big(\theta_{i,j}^m {(s_i^k)}^2 + \varphi_{i,j}^m s_i^k + \eta_{i,j}^m\big) - \mathbb{M}(1-x_{i,j}^{m,k}), \nonumber\\
\IEEEeqnarraymulticol{3}{c}{\forall\, i \in \mathcal{V}^-,~ j \in \mathcal{V}^+,~ m \in T_{i,j}.} \label{eq:7}
\end{IEEEeqnarray}

\begin{IEEEeqnarray}{rCl}
s_i^k &=& \sum_{j \in \mathcal{V}-\{0,i\}} \sum_{m \in T_{i,j}} w_{i,j}^{m,k}, \nonumber\\
\IEEEeqnarraymulticol{3}{c}{\forall i \in \mathcal{V}^-.} 
\end{IEEEeqnarray}

\begin{IEEEeqnarray}{rCl}
t_{i,j}^{m-1} x_{i,j}^{m,k} &\le& w_{i,j}^{m,k} \le t_{i,j}^{m} x_{i,j}^{m,k}, \nonumber\\
\IEEEeqnarraymulticol{3}{c}{\forall (i,j) \in \mathcal{B},~ m \in T_{i,j}.} 
\end{IEEEeqnarray}

\end{small}

\section{Meta-heuristic method}

\begin{figure}[t]
\centering
\begin{minipage}{0.45\textwidth}  
\includegraphics[width=\textwidth]{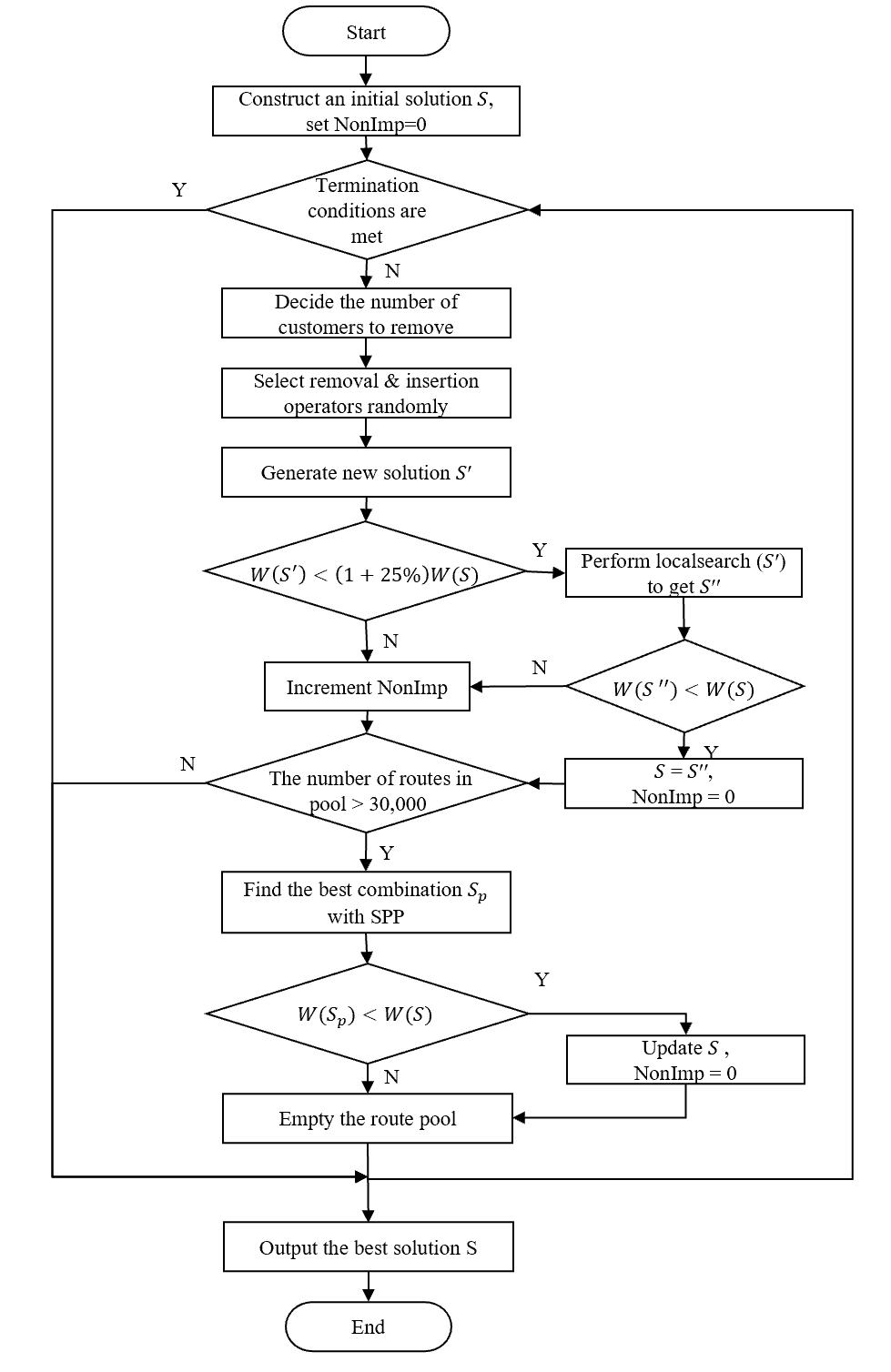}  
\end{minipage}
\caption{The entire scheme of our LNS-SPP method}
\label{fig:LNS-SPP}
\end{figure}

This section describes a hybrid metaheuristic designed for the ALD-EVRP, which combines a Large Neighborhood Search (LNS) framework, a Local Search (LS) procedure, and a Set Partitioning Problem (SPP) component \citep{zhao2022adaptive}. The LNS algorithm iteratively adjusts the solution using multiple removal and insertion operators, while the LS procedure further improves the solution by exploring its neighborhood. The routes generated from LNS and LS are stored in a route pool, and once its size surpasses a predefined limit, the SPP module is triggered to identify the optimal combination of routes that visits all customers exactly once. The overall framework, illustrated in Figure~\ref{fig:LNS-SPP}.

The LNS component dynamically explores diverse neighborhoods by iteratively removing and reinserting customer nodes, enabling effective restructuring of routes to reduce total energy consumption. Our framework incorporates a suite of customized removal operators designed specifically for the characteristics of the ALD-EVRP. The distance-based removal operator identifies nodes that are farthest from their current tours, effectively targeting outliers that may cause inefficiencies. The load-based removal prioritizes nodes with heavier loads or later service times, which helps balance vehicle utilization and mitigate energy peaks. The Shaw removal, in contrast, removes clusters of nodes with similar attributes, such as proximity or demand level, to introduce greater diversity and escape local optima. Following the removal phase, a set of insertion operators—including Random, Greedy, and Regret-based strategies—reconstructs feasible solutions by sequentially reinserting nodes in a way that balances exploration of the search space with solution quality.

To further refine solutions, the local search (LS) procedure intensively examines promising candidates within their neighborhoods, employing problem-specific operators such as node swaps, relocations, and segment reorganization. These operators exploit structural properties of the routes, for example, by reordering delivery sequences to reduce detours or smoothing speed profiles to lower energy consumption. All routes generated during the LNS and LS phases are maintained in a route pool, which serves as a repository of high-quality sub-solutions. When the pool exceeds a predefined threshold, the SPP module integrates selected sub-routes to ensure global consistency across the solution set while controlling computational complexity. This integration allows the algorithm to combine strong partial solutions from different neighborhoods, effectively leveraging accumulated knowledge from the search process. Overall, this framework provides a flexible yet powerful approach for efficiently optimizing large-scale ALD-EVRP instances under realistic operational constraints.


\section{Computational Experiment}
This section introduces a set of newly generated benchmark instances for the ALD-EVRP and reports computational experiments designed to assess the effectiveness of the proposed LNS-SPP algorithm. The algorithm is coded in Python, while small-scale instances are solved using the Bonmin solver.

\subsection{Experiments on Small Scale Instances:} 
We generated small scale datasets by sub-sampling the original data, consisting of 4, 7, 10, and 15 customers with a single vehicle. These datasets were solved through the Bonmin solver and our heuristic method, LNS-LS. Although BonMin is capable of solving general Mixed Integer NonLinear Programming problems, the computational difficulty arises from the high-dimensional product of three variables \(W\), which imposes significant challenges for exact solvers. In applying BonMin, we disregarded the effect of load on energy consumption.

The results are summarized in Table~\ref{tab:BonMin}. LNS-LS consistently achieves the same optimal solutions as Bonmin or identifies improved feasible solutions for instances that Bonmin cannot solve, while requiring substantially less runtime. This demonstrates the efficiency of our algorithm. Furthermore, the algorithm can handle more complex objective functions, scales well with problem size, and can be generalized to broader EVRP variants with diverse constraints.
\begin{table}[ht]
\caption{Comparison of solutions obtained by LNS-LS and Bonmin on small instances}
\label{tab:BonMin}
\centering
\resizebox{0.5\textwidth}{!}{
\begin{tabular}{rrrrrrrrrr} 
\toprule
\multirow{2}{*}{n} & \multicolumn{3}{c}{\textbf{BonMin}} & \multicolumn{1}{c}{} & \multicolumn{5}{c}{\textbf{LNS-LS}} \\ 
\cline{2-4}\cline{6-10}
 & \multicolumn{1}{c}{O/F} & \multicolumn{1}{c}{$W$} & \multicolumn{1}{c}{$T(s)$} & \multicolumn{1}{c}{} & \multicolumn{1}{c}{$W_{best}$} & \multicolumn{1}{c}{Gap} & \multicolumn{1}{c}{TT} & \multicolumn{1}{c}{$W_{avg}$} & \multicolumn{1}{c}{$T(s)$} \\ 
\hline
\textbf{4A} & O & 63275.08 & 186.08 &  & 63275.08 & 0.00 & 90.91 & 63275.08 & 1.7 \\
\textbf{4B} & O & 58150.03 & 947.19 &  & 58150.03 & 0.00 & 82.48 & 58150.03 & 1.9 \\
\textbf{4C} & O & 65110.54 & 1012.56 &  & 65110.54 & 0.00 & 92.73 & 65110.54 & 0.9 \\
\textbf{4D} & O & 58741.90 & 20500 &  & 58741.90 & 0.00 & 82.60 & 58741.90 & 1.3 \\
\textbf{4E} & O & 61698.88 & 6851.18 &  & 61698.88 & 0.00 & 87.94 & 61698.88 & 1.6 \\ 
\textbf{7A} & O & 65114.77 & 11026.66 &  & 65114.77 & 0.00 & 94.43 & 65114.77 & 1.3 \\
\textbf{7B} & O & 63160.91 & 9911.15 &  & 63160.91 & 0.00 & 91.15 & 63160.91 & 2.0 \\
\textbf{7C} & O & 65410.52 & 10404.99 &  & 65410.52 & 0.00 & 93.47 & 65410.52 & 0.9 \\
\textbf{7D} & O & 59816.67 & 36000 &  & 59816.67 & 0.00 & 85.35 & ~59816.67 & 1.7 \\
\textbf{7E} & O & 65258.60 & 11026.24 &  & 65258.60 & 0.00 & 91.24 & 65258.60 & 1.4 \\ 
\textbf{10A} & O & 66608.74 & 35566.91 &  & 66608.74 & 0.00 & 96.25 & 66608.74 & 1.5 \\
\textbf{10B} & F & 63356.68 & 43200 &  & \textbf{63299.66} & \textbf{-0.09} & 90.41 & 63299.66 & 2.1 \\
\textbf{10C} & O & 66717.67 & 37192.46 &  & 66717.67 & 0.00 & 95.91 & 66717.67 & 1.1 \\
\textbf{10D} & F & 63834.21 & 43200 &  & \textbf{60591.43} & \textbf{-5.08} & 85.88 & 60591.43 & 2.2 \\
\textbf{10E} & F & 65808.75 & 43200 &  & \textbf{65742.94} & \textbf{-0.10} & 92.53 & ~65742.94 & 2.3 \\ 
\textbf{15A} & F & 69088.61 & 43200 &  & \textbf{68363.18} & \textbf{-1.05} & 98.74 & 68363.18 & 2.5 \\
\textbf{15B} & F & 63703.63 & 43200 &  & \textbf{63678.15} & \textbf{-0.04} & 90.55 & 63679.26 & 2.3 \\
\textbf{15C} & F & 76361.37 & 43200 &  & \textbf{68022.71} & \textbf{-10.92} & 96.39 & 68022.71 & 2.8 \\
\textbf{15D} & F & 66628.98 & 43200 &  & \textbf{66282.51} & \textbf{-0.52} & 95.67 & 66285.30 & 3.1 \\
\textbf{15E} & F & 77065.52 & 43200 &  & \textbf{68264.64} & \textbf{-11.42} & 98.14 & 68276.39 & 2.8 \\
\bottomrule
\end{tabular}
}

\end{table}

\subsection{Experiments on Large Scale Instances}
We conduct further experiments on $47$ ALD-EVRP instances, where customer counts range from 30 to 114, to assess the performance of our algorithm.
To generate the two-dimensional distance matrix, we first defined six time intervals based on the provided speed settings, where the vehicle maintains a constant speed in each interval. For example, the vehicle travels at 1.1 km/min between 30 and 60 minutes, and at 0.9 km/min between 90 and 150 minutes. Using these settings, we calculated the speed for 100 evenly spaced time points. To compute the distances, we applied a cumulative approach: for each pair of time points, the time difference between consecutive points is multiplied by the corresponding speed in that interval to accumulate the total distance traveled. This ensures that the resulting matrix adheres to the FIFO  principle, meaning that vehicles departing earlier will always arrive sooner, maintaining the intended sequence.
we solve the real-world instances using the LNS-SPP algorithm with and without considering the real-time load in the EC function, as well as the one with static loads (the initial load). These results are denoted as $W$, $W_{\mathrm{no-load}}$, and $W_{\mathrm{ini-load}}$ respectively in Table~\ref{tab:load}.
\begin{table*}[t]
\caption{Evaluation of load-dependent condition.}
\label{tab:load}
\centering
\resizebox{1\textwidth}{!}{
\begin{tabular}{rrrrrrr|lrrrrrr} 
\toprule
\multicolumn{1}{c}{} & \multicolumn{1}{c}{} & \multicolumn{2}{c}{Evaluation process} & \multicolumn{2}{c}{GAP (result)} & \multicolumn{1}{c|}{} &  & \multicolumn{1}{c}{} & \multicolumn{1}{c}{} & \multicolumn{2}{c}{Evaluation process} & \multicolumn{2}{c}{GAP (result)} \\
\multicolumn{1}{c}{\textbf{Ins}} & \multicolumn{1}{c}{\textbf{n}} & \multicolumn{1}{c}{$W$} & \multicolumn{1}{c}{$W_{no-load}$} & \multicolumn{1}{c}{$G_{no-load}$} & \multicolumn{1}{c}{$G_{ini-load}$} & \multicolumn{1}{c|}{} &  & \multicolumn{1}{c}{\textbf{Ins}} & \multicolumn{1}{c}{\textbf{n}} & \multicolumn{1}{c}{$W$} & \multicolumn{1}{c}{$W_{no-load}$} & \multicolumn{1}{c}{$G_{no-load}$} & \multicolumn{1}{c}{$G_{ini-load}$} \\ 
\cline{1-6}\cline{8-14}
\textbf{001} & 30 & 517857.06 & 506183.39 & -2.25 & \textbf{0.33} &  &  & \textbf{032} & 87 & 1205223.46 & 1179416.64 & -2.14 & \textbf{4.13} \\
\textbf{002} & 41 & 696986.63 & 685460.40 & -1.65 & \textbf{0.14} &  &  & \textbf{033} & 87 & 1053836.57 & 1043644.10 & -0.97 & \textbf{4.34} \\
\textbf{003} & 43 & 524189.38 & 513378.37 & -2.06 & \textbf{0.64} &  &  & \textbf{034} & 87 & 1399988.53 & 1371021.22 & -2.07 & \textbf{4.97} \\
\textbf{004} & 44 & 766585.00 & 751516.05 & -1.97 & \textbf{0.49} &  &  & \textbf{035} & 88 & 1068404.97 & 1040126.44 & -2.65 & \textbf{5.11} \\
\textbf{005} & 51 & 1001786.61 & 980194.27 & -2.16 & \textbf{0.52} &  &  & \textbf{036} & 90 & 1276698.42 & 1235350.49 & -3.24 & \textbf{5.13} \\
\textbf{006} & 53 & 633356.04 & 622590.37 & -1.70 & \textbf{0.47} &  &  & \textbf{037} & 91 & 1206481.57 & 1169729.37 & -3.05 & \textbf{5.41} \\
\textbf{007} & 54 & 496666.94 & 486780.96 & -1.99 & \textbf{1.25} &  &  & \textbf{038} & 93 & 941496.963 & 921385.35 & -2.14 & \textbf{5.19} \\
\textbf{008} & 57 & 879079.49 & 862010.64 & -1.94 & \textbf{1.59} &  &  & \textbf{039} & 98 & 1531762.13 & 1509712.47 & -1.44 & \textbf{5.82} \\
\textbf{009} & 59 & 639686.61 & 633324.97 & -0.99 & \textbf{1.15} &  &  & \textbf{040} & 98 & 1692423.57 & 1642543.95 & -2.95 & \textbf{6.45} \\
\textbf{010} & 60 & 1067694.06 & 1049034.20 & -1.75 & \textbf{1.69} &  &  & \textbf{041} & 104 & 1120009.66 & 1050908.34 & -6.17 & \textbf{2.54} \\
\textbf{011} & 61 & 937753.26 & 918099.49 & -2.10 & \textbf{1.64} &  &  & \textbf{042} & 105 & 1400651.05 & 1327930.50 & -5.19 & \textbf{6.13} \\
\textbf{012} & 62 & 693719.93 & 682636.95 & -1.60 & \textbf{2.34} &  &  & \textbf{043} & 106 & 1416699.78 & 1385438.23 & -2.21 & \textbf{4.59} \\
\textbf{013} & 64 & 990999.71 & 974493.16 & -1.67 & \textbf{1.84} &  &  & \textbf{044} & 110 & 1021367.59 & 990055.39 & -3.07 & \textbf{5.54} \\
\textbf{014} & 64 & 1082607.11 & 1061854.60 & -1.92 & \textbf{2.46} &  &  & \textbf{045} & 110 & 1178598.16 & 1168171.03 & -0.88 & 0.00 \\
\textbf{015} & 65 & 937852.26 & 920341.99 & -1.87 & \textbf{2.40} &  &  & \textbf{046} & 110 & 1534117.14 & 1416467.17 & -7.67 & \textbf{6.07} \\
\textbf{016} & 66 & 784678.40 & 768616.51 & -2.05 & 0.00 &  &  & \textbf{047} & 114 & 1306986.80 & 1268430.06 & -2.95 & \textbf{5.97} \\
\bottomrule
\end{tabular}
}
\end{table*}
The percentage gap $G_{\mathrm{no-load}}$ between $W_{\mathrm{no-load}}$ and $W$ is calculated as $G_{\mathrm{no-load}} = (\frac{W_{\mathrm{no-load}}}
{W}-1)*100\%$, which is always negative, reflecting the over-optimistic estimation of EC when load is ignored.

To ensure fair comparison, we re-evaluated the optimal routes obtained under the static-load and no-load settings by recalculating their EC values under the real-time load model. We then computed the relative gaps $G_{\mathrm{ini}}$ with respect to $W$. Results show that $G_{\mathrm{ini-load}}$ is strictly positive in 30 out of 32 instances, confirming that  the static-load model tends to overestimate EC, as the actual load decreases along the trip. These results highlight the importance of explicitly incorporating real-time load into the EC function.

As reported in Table~\ref{tab:load}, LNS-SPP with real-time load consideration ($W$) consistently produces reliable and high-quality solutions across all large-scale instances, scaling up to more than 100 customers. Furthermore, the algorithm successfully solves all instances within reasonable computation time, whereas exact solvers such as BonMin fail to handle problems of this size.

Overall, the results show that the proposed LNS-SPP scales well to large ALD-EVRP instances and outperforms baseline models, emphasizing the benefit of considering real-time load in the energy consumption function.

\section{Conclusion}
In this paper, we present a new time-dependent routing problem for electric vehicles (ALD-EVRP) that incorporates real-time acceleration and load dynamics. We introduce a realistic speed model that considers acceleration effects, in contrast to the standard stepwise speed function. The mathematical model provides a more accurate representation of real-world electric vehicle operations.

To address the computational challenges of this complex problem, we develop the LNS-SPP algorithm that effectively handles the ALD-EVRP. The algorithm demonstrates excellent scalability, solving instances with up to 114 customers while maintaining high solution quality. 

Future work will focus on extending the model to handle multiple vehicle types and dynamic traffic conditions, further enhancing the practical applicability of our approach.

\bibliographystyle{IEEEtran} 
\bibliography{aldvrp} 

@article{KANCHARLA2020113714,
title = {Electric vehicle routing problem with non-linear charging and load-dependent discharging},
journal = {Expert Systems with Applications},
volume = {160},
pages = {113714},
year = {2020},
issn = {0957-4174},
doi = {https://doi.org/10.1016/j.eswa.2020.113714},
url = {https://www.sciencedirect.com/science/article/pii/S0957417420305388},
author = {Surendra Reddy Kancharla and Gitakrishnan Ramadurai}
}

@article{BEHNKE2021794,
title = {A column generation approach for an emission-oriented vehicle routing problem on a multigraph},
journal = {European Journal of Operational Research},
volume = {288},
number = {3},
pages = {794-809},
year = {2021},
issn = {0377-2217},
doi = {https://doi.org/10.1016/j.ejor.2020.06.035},
url = {https://www.sciencedirect.com/science/article/pii/S0377221720305786},
author = {Martin Behnke and Thomas Kirschstein and Christian Bierwirth}
}

@article{GALVIN2017234,
title = {Energy consumption effects of speed and acceleration in electric vehicles: Laboratory case studies and implications for drivers and policymakers},
journal = {Transportation Research Part D: Transport and Environment},
volume = {53},
pages = {234-248},
year = {2017},
issn = {1361-9209},
doi = {https://doi.org/10.1016/j.trd.2017.04.020},
url = {https://www.sciencedirect.com/science/article/pii/S1361920915301280},
author = {Ray Galvin}
}

@article{PRINS20041985,
title = {A simple and effective evolutionary algorithm for the vehicle routing problem},
journal = {Computers \& Operations Research},
volume = {31},
number = {12},
pages = {1985-2002},
year = {2004},
issn = {0305-0548},
doi = {https://doi.org/10.1016/S0305-0548(03)00158-8},
url = {https://www.sciencedirect.com/science/article/pii/S0305054803001588},
author = {Christian Prins}
}

@article{zhao2022adaptive,
  title={Adaptive large neighborhood search for the time-dependent profitable dial-a-ride problem},
  author={Zhao, Jingyi and Poon, Mark and Zhang, Zhenzhen and Gu, Ruixue},
  journal={Computers \& Operations Research},
  volume={147},
  pages={105938},
  year={2022},
  publisher={Elsevier}
}

@article{naeem2024energy,
  author={Naeem, H. M. Y. and Bhatti, A. I. and Butt, Y. A. and Ahmed, Q.},
  title={Energy Economization Using Connectivity-Based Eco-Routing and Driving for Fleet of Battery Electric Vehicles},
  journal={IEEE Transactions on Transportation Electrification},
  year={2024},
  volume={10},
  number={1},
  pages={1923--1934}
}

\clearpage
\begingroup
\onecolumn

\endgroup
\clearpage
\twocolumn

\end{document}